\documentclass{birkjour}
\usepackage{cite}
\usepackage[colorlinks,linkcolor=blue,citecolor=green]{hyperref}
\usepackage{cleveref}

\crefname{equation}{}{}
\newtheorem{theorem}{Theorem}[section]
\newtheorem{lemma}[theorem]{Lemma}

\newtheorem{corollary}[theorem]{Corollary}

\theoremstyle{definition}
\newtheorem{definition}[theorem]{Definition}

\newtheorem{notation}[theorem]{Notation}

\theoremstyle{remark}
\newtheorem{remark}[theorem]{Remark}

\numberwithin{equation}{section}

\begin{document}

%
%
%
%
%
%
%
%
%

\title[Boundary H\"{o}lder Regularity on Reifenberg Flat Domains]
 {Boundary H\"{o}lder Regularity for Elliptic Equations on Reifenberg Flat Domains}

\author{Yuanyuan Lian}

\address{%
School of Mathematical Sciences\\
Shanghai Jiao Tong University\\
Shanghai, 200240\\
PR China}

\email{lianyuanyuan@sjtu.edu.cn; lianyuanyuan.hthk@gmail.com}

\thanks{This research is supported by the China Postdoctoral Science Foundation (Grant No.
2021M692086), the National Natural Science Foundation of China (Grant No. 12031012 and 11831003) and  the Institute of Modern Analysis-A Frontier Research Center of Shanghai.}
\author{Kai Zhang}
\address{School of Mathematical Sciences\\
Shanghai Jiao Tong University\\
Shanghai, 200240\\
PR China}
\email{zhangkaizfz@gmail.com}


\date{May 19, 2021}
\dedicatory{}

\begin{abstract}
In this paper, we investigate the boundary H\"{o}lder regularity for elliptic equations (precisely, the Poisson equation, linear equations in divergence form and non-divergence form, the p-Laplace equations and fully nonlinear elliptic equations) on Reifenberg flat domains. We prove that for any $0<\alpha<1$, there exists $\delta>0$ such that the solution is $C^{\alpha}$ at $x_0\in \partial \Omega$ provided that $\Omega$ is $\delta$-Reifenberg flat at $x_0$ (see \Cref{d-re}). In particular, for any $0 < \alpha < 1$, if $\partial \Omega$ is $C^1$ and $u=g$ on $\partial \Omega$ with $g\in C^{\alpha}(x_0)$, then $u\in C^{\alpha}(x_0)$. A similar result for the Poisson equation has been proved by Lemenant and Sire \cite{MR3156895}, where the Alt-Caffarelli-Friedman's monotonicity formula is used.
\end{abstract}

\maketitle
\section{Introduction}
\label{intro}
In the regularity theory of partial differential equations, the H\"{o}lder continuity is a kind of quantitative estimate. It is usually the first smooth regularity for solutions and the beginning for higher regularity. Take the uniformly elliptic equations in divergence form or non-divergence form for example. With respect to the interior H\"{o}lder regularity, De Giorgi \cite{MR0093649}, Nash \cite{MR0100158} and Moser \cite{MR0170091} proved it for elliptic equations in divergence form, and it was extended to elliptic equations in nondivergence form by Krylov and Safonov \cite{MR525227}. In particular, the fully nonlinear elliptic equations can also be treated (see \cite{MR1005611} and \cite{MR1351007}). With regard to the boundary H\"{o}lder continuity, a well-known result is that if $\Omega$ satisfies the exterior cone condition at $x_0\in \partial \Omega $, the solution is H\"{o}lder continuous at $x_0$ (see \cite{MR0221087} and \cite[Theorem 8.29 and Corollary 9.28]{MR1814364}).

All results above mentioned show that there exists some $0<\alpha<1$ such that $u\in C^{\alpha}$. It is interesting to ask whether the exponent $\alpha$ could be bigger? Based only on the uniform ellipticity condition, it is impossible to obtain a higher interior H\"{o}lder regularity (see \cite{MR882838}). On the other hand, for the boundary H\"{o}lder regularity, we can do expect a bigger exponent.



For the Poisson equation, Lemenant and Sire \cite{MR3156895} proved that for any $0<\alpha<1$, there exists $\delta>0$ such that the solution is $C^{\alpha}$ at $x_0\in \partial \Omega$ provided that $\Omega$ is $\delta$-Reifenberg flat at $x_0$. In this paper, we prove analogue results for different types of elliptic equations and our approach is simple. Moreover, even for the Poisson equation, we relax the requirements imposed in \cite{MR3156895}. The main idea in our method is that a Reifenberg flat domain can be regarded as a perturbation of half balls in different scales.

First, we introduce some notions.
\begin{definition}[\textbf{Reifenberg flat domain}]\label{d-re}
We say that $\Omega$ is $\delta$-Reifenberg flat from the exterior at $x_0\in \partial \Omega$ if there exists $r_0>0$ such that the following holds: for any $0<r<r_0$, there exists a coordinate system $\{y_1,...,y_n \}$ (isometric to the original coordinate system) such that $x_0=0$ in this coordinate system and
\begin{equation}\label{e-re}
B_r\cap \Omega \subset B_r \cap\{y_n >-\delta r\}.
\end{equation}
\end{definition}

Next, we introduce a pointwise characterization for a function.
\begin{definition}\label{d-f}
Let $\Omega \subset R^{n}$ be a bounded set (may be not a domain) and $f:\Omega\rightarrow R$ be a function. We say that $f$ is $C^{\alpha}$ ($0<\alpha\leq 1$) at $x_0\in \Omega$ or $f\in C^{\alpha}(x_0)$ if there exist $r_0>0$ and $K>0$ such that
\begin{equation}\label{holder}
  |f(x)-f(x_0)|\leq K|x-x_0|^{\alpha},~~\forall~x\in \Omega\cap B_{r_0}(x_0).
\end{equation}
Then, define
\begin{equation*}
[f]_{C^{\alpha}(x_0)}=\min \left\{K\big | \cref{holder} ~\mbox{holds with}~K\right\}
\end{equation*}
and
\begin{equation*}
\|f\|_{C^{\alpha}(x_0)}=|f(x_0)|+[f]_{C^{\alpha}(x_0)}.
\end{equation*}

We say that $f$ is $C_{q}^{-p,\alpha}$ at $x_0$ or $f\in C_q^{-p,\alpha}(x_0)$ ($p>0,q\geq 1,\alpha>0$) if there exist $r_0>0$ and $K>0$ such that
\begin{equation}\label{e.c-1}
\|f\|_{L^q(\Omega\cap B_r(x_0) )}\leq Kr^{\alpha-p+n/q}, ~\forall ~0<r<r_0.
\end{equation}
Similarly, we set
\begin{equation*}
\|f\|_{C_{q}^{-p,\alpha}(x_0)}=\min \left\{K \big | \cref{e.c-1} ~\mbox{holds with}~K\right\}.
\end{equation*}
If $f\in C_{q}^{-p, \alpha}(x)$ for any $x\in \Omega$ with the same $r_0$ and
\begin{equation*}
  \|f\|_{C_{q}^{-p,\alpha}(\Omega)}:= \sup_{x\in \Omega} \|f\|_{C_{q}^{-p,\alpha}(x)}<+\infty,
\end{equation*}
we say that $f\in C_{q}^{-p,\alpha}(\Omega)$.
\end{definition}

\begin{remark}\label{r-1}
Since we study the boundary pointwise regularity, throughout this paper, we always assume that $0\in \partial \Omega$ and $r_0=1$ in \Cref{d-re} and \Cref{d-f} without loss of generality.
\end{remark}

Our main results are listed in the following. First, we consider the Laplace operator:
\begin{theorem}\label{t-1}
Let $0<\alpha<1$ and $u$ be a weak solution of
\begin{equation}\label{e-poisson}
\left\{\begin{aligned}
&\Delta u=f&& ~~\mbox{in}~~\Omega\cap B_1;\\
&u=g&& ~~\mbox{on}~~\partial \Omega\cap B_1,
\end{aligned}\right.
\end{equation}
where $g\in C^{\alpha}(0)$ and $f\in C_{p}^{-2,\alpha}(0) $ for some $p>n/2$. Suppose that $\Omega$ is $\delta$-Reifenberg flat from the exterior at $0$, where $\delta>0$ depends only on $n$ and $\alpha$.

Then $u$ is $C^{\alpha}$ at $0$ and
\begin{equation*}\label{e-poisson-holder}
  |u(x)-u(0)|\leq C |x|^{\alpha}\left(\|u\|_{L^{\infty }(\Omega\cap B_1)}+\|f\|_{C_{p}^{-2,\alpha}(x_0)}+[g]_{C^{\alpha}(0)}\right), ~~\forall ~x\in \Omega\cap B_1,
\end{equation*}
where $C$ depends only on $n,\alpha$ and $p$.
\end{theorem}

\begin{remark}\label{r-1.3}
If $f\in L^q$ with $q>n/(2-\alpha)$, then $f\in C_{p}^{-2,\alpha}(0) $ for some $p>n/2$. In particular, $f\in L^n$ is enough for any $0<\alpha<1$.
\end{remark}

\begin{remark}\label{r-1.1}
\Cref{t-1} is more general by comparing with the result of \cite{MR3156895}.
\end{remark}

For linear elliptic equations in divergence form, we have
\begin{theorem}\label{t-2}
Let $0<\alpha<1$ and $u$ be a weak solution of
\begin{equation}\label{e.div}
\left\{\begin{aligned}
&(a^{ij}u_i)_j=f&& ~~\mbox{in}~~\Omega\cap B_1;\\
&u=g&& ~~\mbox{on}~~\partial \Omega\cap B_1,
\end{aligned}\right.
\end{equation}
where $a^{ij}$ is uniformly elliptic with ellipticity constants $\lambda$ and $\Lambda$, $g\in C^{\alpha}(0)$ and $f\in C_{p}^{-2,\alpha}(0) $ for some $p>n/2$. Suppose that $|a^{ij} -\delta^{ij}|\leq \delta$ and $\Omega$ is $\delta$-Reifenberg flat from the exterior at $0$, where $\delta>0$ depends only on $n, \lambda, \Lambda$ and $\alpha$.

Then $u$ is $C^{\alpha}$ at $0$ and
\begin{equation*}
  |u(x)-u(0)|\leq C |x|^{\alpha}\left(\|u\|_{L^{\infty }(\Omega\cap B_1)}+\|f\|_{C_{p}^{-2,\alpha}(0)}+[g]_{C^{\alpha}(0)}\right), ~~\forall ~x\in \Omega\cap B_1,
\end{equation*}
where $C$ depends only on $n,\lambda, \Lambda, \alpha$ and $p$.
\end{theorem}

\begin{remark}\label{r-2}
In \cref{e.div}, the Einstein summation convention is used (similarly hereinafter), i.e., repeated indices means summation. The symbol $\delta^{ij}$ denotes the unit matrix. In fact, $\delta^{ij}$ can be replaced by any constant symmetric matrix with eigenvalues lying in $[\lambda, \Lambda]$.
\end{remark}

\begin{remark}\label{r-1.4}
The smallness condition $|a^{ij}-\delta^{ij}|\leq \delta$ is necessary, which is different from the boundary regularity for equations in non-divergence form (see \Cref{t-3} below). In fact, based only on the uniform ellipticity, we can't expect any higher boundar H\"{o}lder regularity for elliptic equations in divergence form (see \cite[(1.4)]{MR2069724} and \cite[Section 5]{MR159110}).
\end{remark}

For linear equations in non-divergence form, the corresponding boundary H\"{o}lder regularity holds. More generally, it holds for fully nonlinear equations in non-divergence form:
\begin{theorem}\label{t-3}
Let $0<\alpha<1$ and $u$ be a viscosity solution of
\begin{equation}\label{e.ful}
\left\{\begin{aligned}
&u\in S(\lambda,\Lambda,f)&& ~~\mbox{in}~~\Omega\cap B_1;\\
&u=g&& ~~\mbox{on}~~\partial \Omega\cap B_1,
\end{aligned}\right.
\end{equation}
where $g\in C^{\alpha}(0)$ and $f\in L^n(\Omega\cap B_1)$.  Suppose that $\Omega$ is $\delta$-Reifenberg flat from the exterior at $0$, where $\delta>0$ depends only on $n, \lambda, \Lambda$ and $\alpha$.

Then $u$ is $C^{\alpha}$ at $0$ and
\begin{equation}\label{e-holder}
  |u(x)-u(0)|\leq C |x|^{\alpha}\left(\|u\|_{L^{\infty }(\Omega\cap B_1)}
  +\|f\|_{L^{n}(\Omega\cap B_1)}+[g]_{C^{\alpha}(0)}\right), ~~\forall ~x\in \Omega\cap B_1,
\end{equation}
where $C$ depends only on $n, \lambda, \Lambda$ and $\alpha$.
\end{theorem}

\begin{remark}\label{r-1.2}
In \cref{e.ful}, $S(\lambda,\Lambda,f)$ denotes the Pucci class with ellipticity constants $\lambda, \Lambda$ and righthand $f$. The Pucci class is a natural generalization of linear uniformly elliptic equations in non-divergence form. For the notion of viscosity solutions and basic properties of the Pucci class, we refer to \cite{MR1351007}, \cite{MR1376656} and \cite{MR1118699}.
\end{remark}

\begin{remark}\label{r-1.5}
We can relax $f\in L^n$ to $f\in L^{n-\varepsilon_0}$ for some $\varepsilon_0>0$ depending only on $n,\lambda, \Lambda$ and $\alpha$ (see \cite{MR1237053}) since the proof mainly relies on the A-B-P maximum principle.
\end{remark}

For the $p$-Laplace equations, we have
\begin{theorem}\label{t-4}
Let $0<\alpha<1$, $1<p<\infty$ and $u$ be a weak solution of
\begin{equation}\label{e.p-Lap}
\left\{\begin{aligned}
&\mathrm{div}(|Du|^{p-2}Du)=0&& ~~\mbox{in}~~\Omega\cap B_1;\\
&u=g&& ~~\mbox{on}~~\partial \Omega\cap B_1,
\end{aligned}\right.
\end{equation}
where $g\in C^{\alpha}(0)$. Suppose that $\Omega$ is $\delta$-Reifenberg flat from the exterior at $0$, where $\delta>0$ depends only on $n, \alpha$ and $p$.

Then $u$ is $C^{\alpha}$ at $0$ and
\begin{equation*}
  |u(x)-u(0)|\leq C |x|^{\alpha}\left(\|u\|_{L^{\infty }(\Omega\cap B_1)}+[g]_{C^{\alpha}(0)}\right), ~~\forall ~x\in \Omega\cap B_1,
\end{equation*}
where $C$ depends only on $n, \alpha$ and $p$.
\end{theorem}

For $2\leq p<\infty$, we have the boundary H\"{o}lder regularity for the Poisson equations:
\begin{theorem}\label{t-5}
Let $0<\alpha<1$, $2\leq p<\infty$ and $u$ be a weak solution of
\begin{equation}\label{e.p-Poi}
\left\{\begin{aligned}
&\mathrm{div}(|Du|^{p-2}Du)=f&& ~~\mbox{in}~~\Omega\cap B_1;\\
&u=g&& ~~\mbox{on}~~\partial \Omega\cap B_1,
\end{aligned}\right.
\end{equation}
where $g\in C^{\alpha}(0)$ and $f\in C_{q}^{-p,(p-1)\alpha}(0) $ for some $q>n/p$ and $1/p+1/q\leq 1$. Suppose that $\Omega$ is $\delta$-Reifenberg flat from the exterior at $0$, where $\delta>0$ depends only on $n,\alpha$ and $p$.

Then $u$ is $C^{\alpha}$ at $0$ and
\begin{equation*}
  |u(x)-u(0)|\leq C |x|^{\alpha}\left(\|u\|_{L^{\infty }(\Omega\cap B_1)}+\|f\|_{C_{q}^{-p,(p-1)\alpha}(0)}+[g]_{C^{\alpha}(0)}\right), ~~\forall ~x\in \Omega\cap B_1,
\end{equation*}
where $C$ depends only on $n, \alpha, p$ and $q$.
\end{theorem}

We say that $\Omega$ satisfies the exterior cone condition with slope $\delta$ at $x_0\in \partial \Omega$ if there exist $r_0>0$ and a unit vector $\vec{n}$ such that
\begin{equation*}
 \left\{x\in B(x_0,r): (x-x_0)\cdot \vec{n}<-\delta |x-x_0|\right\}\subset \Omega^c.
\end{equation*}
Clearly, if $\Omega$ satisfies the exterior cone condition with slope $\delta$, it is $\delta$-Reifenberg flat from the exterior. Hence, we have
\begin{corollary}\label{co-1}
In \Crefrange{t-1}{t-5}, the condition that $\Omega$ is $\delta$-Reifenberg flat from the exterior at $0$ can be replaced by that $\Omega$ satisfies the exterior cone condition at $0$ with slope $\delta$. In particular, if $\partial \Omega\in C^1$, \Crefrange{t-1}{t-5} hold.
\end{corollary}

\begin{notation}\label{no1.1}
\begin{enumerate}~~\\
\item $\{e_i\}^{n}_{i=1}$: the standard basis of $R^n$, i.e., $e_i=(0,...0,\underset{i^{th}}{1},0,...0)$.
\item $x'=(x_1,x_2,...,x_{n-1})$ and $x=(x_1,...,x_n)=(x',x_n)$ .
\item $S^n$: the set of $n\times n$ symmetric matrices and $\|A\|=$ the spectral radius of $A$ for any $A\in S^n$.
\item $R^n_+=\{x\in R^n\big|x_n>0\}$.
\item $B_r(x_0)=\{x\in R^{n}\big| |x-x_0|<r\}$, $B_r=B_r(0)$, $B_r^+(x_0)=B_r(x_0)\cap R^n_+$ and $B_r^+=B^+_r(0)$.
\item $T_r(x_0)\ =\{(x',0)\in R^{n}\big| |x'-x_0'|<r\}$ \mbox{ and } $T_r=T_r(0)$.
\item $A^c$: the complement of $A$ and $\bar A $: the closure of $A$, $\forall A\subset R^n$.
\item $\Omega_r=\Omega\cap B_r$ and $(\partial\Omega)_r=\partial\Omega\cap B_r$.
\item $\varphi _i=D_i \varphi=\partial \varphi/\partial x _{i}$ and $D \varphi=(\varphi_1 ,...,\varphi_{n} )$. Similarly, $\varphi _{ij}=D_{ij}\varphi =\partial ^{2}\varphi/\partial x_{i}\partial x_{j}$ and $D^2 \varphi =\left(\varphi _{ij}\right)_{n\times n}$.
\end{enumerate}
\end{notation}

\section{Boundary H\"{o}lder regularity}\label{sec:1}

In this section, we give the detailed proofs of \Crefrange{t-1}{t-5}. The main idea is the following. Since we consider the continuity of the solution up to the boundary, we only need to control the solution from both the above and the below. This is usually carried out by constructing a proper barrier. Here, we construct a sequence of ``barrier-like'' functions and adopt the scaling argument to prove the continuity up to the boundary. Although we can't obtain the H\"{o}lder regularity at once as done by the method of constructing a barrier, the method is more flexible.

We regard the boundary of a Reifenberg flat domain as a perturbation of a hyperplane in different scales. Hence, we only need to construct ``barrier-like'' functions with respect to a flat boundary. Thus, they are easy to construct. In fact, a function like
\begin{equation}\label{barrier}
v_0(x)=\left(\frac{1}{2}\right)^{-\beta}-\left|x
  +\left(\delta+\frac{1}{2}\right)e_n\right|^{-\beta}
\end{equation}
is enough for \Crefrange{t-1}{t-5} by proper rescaling and choosing a sufficient large $\beta$. From above observation, this kind of boundary regularity is the most easily to prove. One can compare it with the boundary H\"{o}lder regularity for ``some'' rather than ``any'' $0<\alpha<1$ (see \cite{MR4163131}), the boundary Lipschitz regularity (see \cite{Safonov2008} and \cite{lian2018boundary}), the boundary differentiability (see \cite{MR2264018,MR2494685,MR3029135}) and the boundary $C^{k,\alpha}$ regularity for $k\geq 1$ (see \cite{MR4088470} and \cite{Lian-Wang-Zhang2020}).
~\\

In the following, we give the detailed proofs of our main results.

\noindent \textbf{Proof of \Cref{t-1}.} Without loss of generality, we assume that $g(0)=0$. Let $M=\|u\|_{L^{\infty }(\Omega\cap B_1)}+\|f\|_{C_{p}^{-2,\alpha}(x_0)}+[g]_{C^{\alpha}(0)}$ and $\Omega_{r}=\Omega\cap B_{r}$. To prove that $u$ is $C^{\alpha}$ at 0, we only need to show the following:

there exist constants $0<\delta<1/4$, $0< \eta < 1$ (depending only on $n$ and $\alpha$) and $\hat{C}\geq 1$ (depending only on $n,\alpha$ and $p$) such that for all $k\geq 0$,

\begin{equation}\label{e-Lap-discrete}
\|u\|_{L^{\infty }(\Omega _{\eta ^{k}})}\leq \hat{C} M \eta^{k\alpha}.
\end{equation}

We prove \cref{e-Lap-discrete} by induction. For $k=0$, it holds clearly. Suppose that it holds for $k$. We need to prove that it holds for $k+1$.

Let $r=\eta ^{k}/2$ and there exists a new coordinate system denoted by $\{x_1,...,x_n \}$ again such that
\begin{equation}\label{e-re-2}
B_r\cap \Omega \subset B_r \cap\{x_n >-\delta r\}.
\end{equation}
Let $\tilde{B}^{+}_{r}=B^{+}_{r}-\delta r e_n $, $\tilde{T}_r=T_r-\delta r e_n$ and $\tilde{\Omega }_{r}=\Omega \cap \tilde{B}^{+}_{r}$ where $e_n=(0,0,...,0,1)$. Then $\Omega _{r/2}\subset \tilde{\Omega }_{r}\subset \Omega_{\eta^k}$.

Take (note that $v_0$ is defined in \cref{barrier})
\begin{equation*}
v(x)=\hat{C}M\eta^{k\alpha}v_0(x/r)
\end{equation*}
and by choosing $\beta$ large enough, $v$ satisfies
\begin{equation*}
\left\{\begin{aligned}
 &\Delta v\leq  0 &&\mbox{in}~~\tilde{B}^{+}_{r}; \\
 &v\geq 0  &&\mbox{in}~~\tilde{B}^{+}_{r};\\
 &v\geq \hat{C} M \eta^{k\alpha}&&\mbox{on}~~\partial \tilde{B}^{+}_{r}\backslash \tilde{T}_{r};\\
  &v(-\delta r e_n)=0.
\end{aligned}
\right.
\end{equation*}
Set $w=u-v$ and we have (note that $v\geq 0$)
\begin{equation*}
    \left\{
    \begin{aligned}
      &\Delta w\geq f &&\mbox{in}~~ \tilde{\Omega }_{r}; \\
      &w = g-v\leq g &&\mbox{on}~~\partial \Omega \cap \tilde{B}^{+}_{r};\\
      &w\leq 0 &&\mbox{on}~~\partial \tilde{B}^{+}_{r}\cap \bar{\Omega}.
    \end{aligned}
    \right.
\end{equation*}

Let
\begin{equation*}
  \delta = \eta.
\end{equation*}
For $v$, it can be calculated directly that
\begin{equation}\label{e-Lap-v}
v\leq C_1\eta\cdot  \hat{C} M \eta^{k\alpha}= C_1\eta ^{1-\alpha }\cdot \hat{C}M\eta ^{(k+1)\alpha}~\mbox{on}~\left\{(x',x_n): x'=0,-\delta r\leq x_n\leq 2\eta r\right\},
\end{equation}
where $C_1$ depends only on $n$. For $w$, by the A-B-P maximum principle, we have for $x\in \tilde{\Omega }_{r}$,
\begin{equation}\label{e-Lap-w}
  \begin{aligned}
w(x) &\leq\|g\|_{L^{\infty }(\partial \Omega \cap \tilde{B}^{+}_{r})}
+Cr^{2-n/p}\|f\|_{L^p(\tilde\Omega_{r})}\\
&\leq M(2r)^{\alpha}+CM(2r)^{\alpha}\\
&=\frac{C_2}{\hat{C}\eta^{\alpha}}\cdot \hat{C}M\eta^{(k+1)\alpha},
  \end{aligned}
\end{equation}
where $C_2$ depends only on $n$ and $p$.

Take $\eta\leq 1/4$ small enough such that
\begin{equation}\label{e.2.3}
  C_1\eta ^{1-\alpha }\leq 1/2.
\end{equation}
Next, take $\hat{C}$ large enough such that
\begin{equation*}
  \frac{C_2}{\hat{C}\eta^{\alpha}}\leq 1/2.
\end{equation*}
Then by combining with \cref{e-Lap-v} and \cref{e-Lap-w}, we have
\begin{equation*}
 u= v+w\leq \hat{C}M\eta ^{(k+1)\alpha}~~~~\mbox{on}~\left\{(x',x_n): x'=0,-\delta r\leq x_n\leq 2\eta r\right\}.
\end{equation*}

By translating $v$ to proper positions (with $v(x', -\delta re_n)=0$ for some $x'\in B_{\eta^{k+1}}$) and similar arguments, we obtain
\begin{equation*}
 \sup_{\Omega_{\eta ^{k+1}}}u\leq \hat{C}M\eta ^{(k+1)\alpha}.
\end{equation*}
The proof for
\begin{equation*}
\inf_{\Omega _{\eta^{k+1}}} u\geq -\hat{C}M\eta ^{(k+1)\alpha}
\end{equation*}
is similar and we omit it here. Therefore,
\begin{equation*}
\|u\|_{L^{\infty}(\Omega _{\eta^{k+1}})}\leq \hat{C}M\eta ^{(k+1)\alpha}.
\end{equation*}
By induction, the proof is completed. \qed~\\

\begin{remark}\label{r-2.1}
From \cref{e.2.3}, we infer an explicit relation between $\delta$ and $\alpha$:
\begin{equation*}
  \delta^{1-\alpha}=c_0
\end{equation*}
for some $0<c_0<1$ depending only on $n$.
\end{remark}
~\\
\noindent \textbf{Proof of \Cref{t-2}.} The proof is almost the same as that of \Cref{t-1} since \cref{e.div} can regarded as a perturbation of the Laplace equation. As before, we assume that $g(0)=0$. Let $M=\|u\|_{L^{\infty }(\Omega\cap B_1)}+\|f\|_{C_{p}^{-2,\alpha}(x_0)}+[g]_{C^{\alpha}(0)}$ and $\Omega_{r}=\Omega\cap B_{r}$. To prove that $u$ is $C^{\alpha}$ at 0, we only need to show the following:

there exist constants $0<\delta<1/4$, $0< \eta < 1$ (depending only on $n,\lambda, \Lambda$ and $\alpha$) and $\hat{C}\geq 1$ (depending only on $n,\lambda, \Lambda,\alpha$ and $p$) such that for all $k\geq 0$,

\begin{equation}\label{e-div-discrete}
\|u\|_{L^{\infty }(\Omega _{\eta ^{k}})}\leq \hat{C} M \eta^{k\alpha}.
\end{equation}

We prove \cref{e-div-discrete} by induction. For $k=0$, it holds clearly. Suppose that it holds for $k$. We need to prove that it holds for $k+1$.

Let $r=\eta ^{k}/2$ and there exists a new coordinate system denoted by $\{x_1,...,x_n \}$ again such that
\begin{equation*}
B_r\cap \Omega \subset B_r \cap\{x_n >-\delta r\}.
\end{equation*}
Let $\tilde{B}^{+}_{r}=B^{+}_{r}-\delta r e_n $, $\tilde{T}_r=T_r-\delta r e_n$ and $\tilde{\Omega }_{r}=\Omega \cap \tilde{B}^{+}_{r}$. Then $\Omega _{r/2}\subset \tilde{\Omega }_{r}\subset \Omega_{\eta^k}$.

Define
\begin{equation*}
v(x)=\hat{C}M\eta^{k\alpha}v_0(x/r)
\end{equation*}
and by choosing $\beta$ large enough, $v$ satisfies
\begin{equation*}
\left\{\begin{aligned}
 &\Delta v\leq  0 &&\mbox{in}~~\tilde{B}^{+}_{r}; \\
 &v\geq 0  &&\mbox{in}~~\tilde{B}^{+}_{r};\\
 &v\geq \hat{C} M \eta^{k\alpha}&&\mbox{on}~~\partial \tilde{B}^{+}_{r}\backslash \tilde{T}_{r};\\
  &v(-\delta r e_n)=0.
\end{aligned}
\right.
\end{equation*}
Set $w=u-v$ and we have
\begin{equation*}
    \left\{
    \begin{aligned}
      &(a^{ij}w_i)_j =f+((\delta^{ij}-a^{ij})v_i)_j &&\mbox{in}~~ \tilde{\Omega }_{r}; \\
      &w\leq g &&\mbox{on}~~\partial \Omega \cap \tilde{B}^{+}_{r};\\
      &w\leq 0 &&\mbox{on}~~\partial \tilde{B}^{+}_{r}\cap \bar{\Omega}.
    \end{aligned}
    \right.
\end{equation*}

Take $ \delta = \eta$. As before, for $v$, it can be calculated directly that
\begin{equation}\label{e-div-v}
v\leq C_1\eta ^{1-\alpha }\cdot \hat{C}M\eta ^{(k+1)\alpha}~\mbox{on}~\left\{(x',x_n): x'=0,-\delta r\leq x_n\leq 2\eta r\right\},
\end{equation}
where $C_1$ depends only on $n$. For $w$, by the A-B-P maximum principle, we have for $x\in \tilde{\Omega }_{r}$,
\begin{equation}\label{e-div-w}
  \begin{aligned}
w(x) &\leq\|g\|_{L^{\infty }(\partial \Omega \cap \tilde{B}^{+}_{r})}
+Cr^{2-n/p}\|f\|_{L^p(\tilde\Omega_{r})}+C\eta r\|Dv\|_{L^{\infty}(\tilde\Omega_{r})}\\
&\leq M(2r)^{\alpha}+CM(2r)^{\alpha}+C\eta\cdot\hat{C}M\eta ^{k\alpha}\\
&=\left(\frac{C_2}{\hat{C}\eta^{\alpha}}+C_3\eta^{1-\alpha}\right)\cdot\hat{C}M\eta^{(k+1)\alpha},
  \end{aligned}
\end{equation}
where $C_2$ depends only on $n,\lambda,\Lambda$ and $p$, and $C_3$ depends only on $n,\lambda$ and $\Lambda$.

Take $\eta\leq 1/4$ small enough such that
\begin{equation*}
  (C_1+C_3)\eta ^{1-\alpha }\leq 1/4.
\end{equation*}
Next, take $\hat{C}$ large enough such that
\begin{equation*}
  \frac{C_2}{\hat{C}\eta^{\alpha}}\leq 1/4.
\end{equation*}
Then by combining with \cref{e-div-v} and \cref{e-div-w}, we have
\begin{equation*}
 u= v+w\leq \hat{C}M\eta ^{(k+1)\alpha}~~~~\mbox{on}~\left\{(x',x_n): x'=0,-\delta r\leq x_n\leq 2\eta r\right\}.
\end{equation*}
By translating $v$ to proper positions and similar arguments, we obtain
\begin{equation*}
 \sup_{\Omega_{\eta ^{k+1}}}u\leq \hat{C}M\eta ^{(k+1)\alpha}.
\end{equation*}
The proof for
\begin{equation*}
\inf_{\Omega _{\eta^{k+1}}} u\geq -\hat{C}M\eta ^{(k+1)\alpha}
\end{equation*}
is similar and we omit here. Therefore,
\begin{equation*}
\|u\|_{L^{\infty}(\Omega _{\eta^{k+1}})}\leq \hat{C}M\eta ^{(k+1)\alpha}.
\end{equation*}
By induction, the proof is completed. \qed~\\

\noindent \textbf{Proof of \Cref{t-3}.} Without loss of generality, we assume that $g(0)=0$. Let $M=\|u\|_{L^{\infty }(\Omega\cap B_1)}+\|f\|_{L^{n}(\Omega\cap B_1)}+[g]_{C^{\alpha}(0)}$ and $\Omega_{r}=\Omega\cap B_{r}$. To prove that $u$ is $C^{\alpha}$ at 0, we only need to show the following:

there exist constants $0<\delta<1/4$, $0< \eta < 1$ and $\hat{C}\geq 1$ depending only on $n, \lambda, \Lambda$ and $\alpha$ such that for all $k\geq 0$,

\begin{equation}\label{e-ful-discrete}
\|u\|_{L^{\infty }(\Omega _{\eta ^{k}})}\leq \hat{C} M \eta^{k\alpha}.
\end{equation}

We prove \cref{e-ful-discrete} by induction. For $k=0$, it holds clearly. Suppose that it holds for $k$. We need to prove that it holds for $k+1$.

Let $r=\eta ^{k}/2$ and there exists a new coordinate system denoted by $\{x_1,...,x_n \}$ again such that
\begin{equation*}
B_r\cap \Omega \subset B_r \cap\{x_n >-\delta r\}.
\end{equation*}
Let $\tilde{B}^{+}_{r}=B^{+}_{r}-\delta r e_n $, $\tilde{T}_r=T_r-\delta r e_n$ and $\tilde{\Omega }_{r}=\Omega \cap \tilde{B}^{+}_{r}$. Then $\Omega _{r/2}\subset \tilde{\Omega }_{r}\subset \Omega_{\eta^k}$.

Define
\begin{equation*}
v(x)=\hat{C}M\eta^{k\alpha}v_0(x/r)
\end{equation*}
and by choosing $\beta$ large enough, $v$ satisfies
\begin{equation*}
\left\{\begin{aligned}
 &M^{+}(D^2v,\lambda,\Lambda)\leq 0 &&\mbox{in}~~\tilde{B}^{+}_{r}; \\
 &v\geq 0  &&\mbox{in}~~\tilde{B}^{+}_{r};\\
 &v\geq \hat{C} M \eta^{k\alpha}&&\mbox{on}~~\partial \tilde{B}^{+}_{r}\backslash \tilde{T}_{r};\\
  &v(-\delta r e_n)=0.
\end{aligned}
\right.
\end{equation*}
Set $w=u-v$ and we have
\begin{equation*}
    \left\{
    \begin{aligned}
      &w\in \underline{S}(\lambda,\Lambda , f) &&\mbox{in}~~ \tilde{\Omega }_{r}; \\
      &w\leq g &&\mbox{on}~~\partial \Omega \cap \tilde{B}^{+}_{r};\\
      &w\leq 0 &&\mbox{on}~~\partial \tilde{B}^{+}_{r}\cap \bar{\Omega}.
    \end{aligned}
    \right.
\end{equation*}

Take $\delta = \eta$. For $v$, it can be calculated directly that
\begin{equation}\label{e-v}
v\leq C_1\eta ^{1-\alpha }\cdot \hat{C}M\eta ^{(k+1)\alpha}~\mbox{on}~\left\{(x',x_n): x'=0,-\delta r\leq x_n\leq 2\eta r\right\},
\end{equation}
where $C_1$ depends only on $n,\lambda$ and $\Lambda$. For $w$, by the A-B-P maximum principle, we have for $x\in \tilde{\Omega }_{r}$,
\begin{equation}\label{e-w}
w(x) \leq\|g\|_{L^{\infty }(\partial \Omega \cap \tilde{B}^{+}_{r})}+Cr\|f\|_{L^n(\tilde\Omega_{r})}
\leq M(2r)^{\alpha}+CMr=\frac{C_2}{\hat{C}\eta^{\alpha}}\cdot \hat{C}M\eta^{(k+1)\alpha},
\end{equation}
where $C_2$ depends only on $n,\lambda$ and $\Lambda$.

Take $\eta\leq 1/4$ small enough such that
\begin{equation*}
  C_1\eta ^{1-\alpha }\leq 1/2.
\end{equation*}
Next, take $\hat{C}$ large enough such that
\begin{equation*}
  \frac{C_2}{\hat{C}\eta^{\alpha}}\leq 1/2.
\end{equation*}
Then by combining with \cref{e-v} and \cref{e-w}, we have
\begin{equation*}
 u= v+w\leq \hat{C}M\eta ^{(k+1)\alpha}~~~~\mbox{on}~\left\{(x',x_n): x'=0,-\delta r\leq x_n\leq 2\eta r\right\}.
\end{equation*}

By translating $v$ to proper positions and similar arguments, we obtain
\begin{equation*}
 \sup_{\Omega_{\eta ^{k+1}}}u\leq \hat{C}M\eta ^{(k+1)\alpha}.
\end{equation*}
The proof for
\begin{equation*}
\inf_{\Omega _{\eta^{k+1}}} u\geq -\hat{C}M\eta ^{(k+1)\alpha}
\end{equation*}
is similar and we omit here. Therefore,
\begin{equation*}
\|u\|_{L^{\infty}(\Omega _{\eta^{k+1}})}\leq \hat{C}M\eta ^{(k+1)\alpha}.
\end{equation*}
By induction, the proof is completed. \qed~\\

For the $p$-Laplace equation, we can prove the boundary H\"{o}lder regularity by the same argument since it can be rewritten as a uniformly elliptic equation in nondivergence form if the gradient of the solution is nonvanishing.
~\\

\noindent \textbf{Proof of \Cref{t-4}.} Without loss of generality, we assume that $g(0)=0$. Let $M=\|u\|_{L^{\infty }(\Omega\cap B_1)}+[g]_{C^{\alpha}(0)}$ and $\Omega_{r}=\Omega\cap B_{r}$. To prove that $u$ is $C^{\alpha}$ at 0, we only need to show the following:

there exist constants $0<\delta<1/4$, $0< \eta < 1$ and $\hat{C}\geq 1$ depending only on $n,\alpha$ and $p$ such that for all $k\geq 0$,

\begin{equation}\label{e-p-Lap-discrete}
\|u\|_{L^{\infty }(\Omega _{\eta ^{k}})}\leq \hat{C} M \eta^{k\alpha}.
\end{equation}

We prove \cref{e-p-Lap-discrete} by induction. For $k=0$, it holds clearly. Suppose that it holds for $k$. We need to prove that it holds for $k+1$.

Let $r=\eta ^{k}/2$. Then there exists a new coordinate system denoted by $\{x_1,...,x_n \}$ again such that
\begin{equation*}
B_r\cap \Omega \subset B_r \cap\{x_n >-\delta r\}.
\end{equation*}
Let $\tilde{B}^{+}_{r}=B^{+}_{r}-\delta r e_n $, $\tilde{T}_r=T_r-\delta r e_n$ and $\tilde{\Omega }_{r}=\Omega \cap \tilde{B}^{+}_{r}$. Then $\Omega _{r/2}\subset \tilde{\Omega }_{r}\subset \Omega_{\eta^k}$.

Define
\begin{equation*}
v(x)=\hat{C}M\eta^{k\alpha}v_0(x/r)
\end{equation*}
and by choosing $\beta$ large enough, $v$ satisfies
\begin{equation*}
\left\{\begin{aligned}
 &\mathrm{div}(|Dv|^{p-2}Dv)\leq 0 &&\mbox{in}~~\tilde{B}^{+}_{r}; \\
 &v\geq 0  &&\mbox{in}~~\tilde{B}^{+}_{r};\\
 &v\geq \hat{C} M \eta^{k\alpha}&&\mbox{on}~~\partial \tilde{B}^{+}_{r}\backslash \tilde{T}_{r};\\
  &v(-\delta r e_n)=0.
\end{aligned}
\right.
\end{equation*}

Take $\delta = \eta$. As before,
\begin{equation}\label{e-p-Lap-v}
v\leq C\eta ^{1-\alpha }\cdot \hat{C}M\eta ^{(k+1)\alpha}~\mbox{on}~\left\{(x',x_n): x'=0,-\delta r\leq x_n\leq 2\eta r\right\},
\end{equation}
where $C$ depends only on $n$ and $p$. From the comparison principle, we have for $x\in \tilde{\Omega }_{r}$,
\begin{equation}\label{e-p-Lap-w}
u(x) \leq v(x)+\|g\|_{L^{\infty }(\partial \Omega \cap \tilde{B}^{+}_{r})}
\leq \left(C\eta ^{1-\alpha }+\frac{C}{\hat{C}\eta^{\alpha}}\right)\cdot \hat{C}M\eta^{(k+1)\alpha},
\end{equation}
where $C$ depends only on $n$ and $p$.

Take $\eta\leq 1/4$ small enough and $\hat{C}$ large enough such that
\begin{equation*}
\left(C\eta ^{1-\alpha }+\frac{C}{\hat{C}\eta^{\alpha}}\right)\leq 1.
\end{equation*}
Then
\begin{equation*}
 u\leq \hat{C}M\eta ^{(k+1)\alpha}~~~~\mbox{on}~\left\{(x',x_n): x'=0,-\delta r\leq x_n\leq 2\eta r\right\}.
\end{equation*}

By translating $v$ to proper positions and similar arguments, we obtain
\begin{equation*}
 \sup_{\Omega_{\eta ^{k+1}}}u\leq \hat{C}M\eta ^{(k+1)\alpha}.
\end{equation*}
The proof for
\begin{equation*}
\inf_{\Omega _{\eta^{k+1}}} u\geq -\hat{C}M\eta ^{(k+1)\alpha}
\end{equation*}
is similar and we omit here. Therefore,
\begin{equation*}
\|u\|_{L^{\infty}(\Omega _{\eta^{k+1}})}\leq \hat{C}M\eta ^{(k+1)\alpha}.
\end{equation*}
By induction, the proof is completed. \qed~\\

For the corresponding Poisson equations, since the difference of two solutions is no longer a solution of some equation, we need to do a little more work. First, we introduce two lemmas, which are motivated by \cite{MR2988770}.
\begin{lemma}\label{l-1}
Let $\phi : [0,\infty)\rightarrow [0,\infty)$ be a nonnegative nonincreasing function. Assume that for some constants $c>0$, $\beta>1$ and $\gamma>0$,
\begin{equation*}
  \phi (h)\leq \frac{c}{(h-k)^{\gamma}} \phi^{\beta}(k), ~\forall ~h>k\geq 0.
\end{equation*}
Then
\begin{equation*}
  \phi(r)=0, ~\mbox{where}~
  r=\left(c2^{\frac{\beta\gamma}{\beta-1}}\phi^{\beta-1}(0)\right)^{1/\gamma}.
\end{equation*}
\end{lemma}
\proof Let $k_m=r(1-2^{-m})$ for $m\geq 0$. Obviously, it is enough to prove that for any $m\geq 0$,
\begin{equation}\label{e.2.1}
\phi (k_m)\leq \frac{\phi (0)}{2^{\mu m}}, ~\mbox{where}~\mu =\frac{\gamma}{\beta-1}.
\end{equation}

We prove \cref{e.2.1} by induction. For $m=0$, it holds clearly. Suppose that it holds for $m$ and we need to prove that it holds for $m+1$. By direct calculation,
\begin{equation*}
  \begin{aligned}
\phi(k_{m+1})&\leq \frac{c}{(2^{-m-1}r)^{\gamma}} \phi^{\beta}(k_m)\\
&\leq 2^{(m+1)\gamma}2^{-\frac{\beta\gamma}{\beta-1}}\phi^{1-\beta}(0)
\left(\frac{\phi(0)}{2^{m\mu}}\right)^{\beta}\\
&= \frac{\phi(0)}{2^{m\mu}}\cdot 2^{(m+1)\gamma-\frac{\beta\gamma}{\beta-1}-m\mu (\beta-1)}\\
&=\frac{\phi(0)}{2^{m\mu}}\cdot 2^{-\frac{\gamma}{\beta-1}}\\
&=\frac{\phi(0)}{2^{(m+1)\mu}}.
  \end{aligned}
\end{equation*}
By induction, the proof is completed.\qed~\\

The next lemma is a kind of A-B-P estimate for the difference of two solutions.
\begin{lemma}\label{l-2}
Let $2\leq p<\infty$. Suppose that $u$ and $v$ are weak solutions of
\begin{equation*}
\left\{\begin{aligned}
&\mathrm{div}(|Du|^{p-2}Du)\geq f_1&& ~~\mbox{in}~~\Omega;\\
&u\leq g_1&& ~~\mbox{on}~~\partial \Omega
\end{aligned}\right.
\end{equation*}
and
\begin{equation*}
\left\{\begin{aligned}
&\mathrm{div}(|Dv|^{p-2}Dv)\leq f_2&& ~~\mbox{in}~~\Omega;\\
&v\geq g_2&& ~~\mbox{on}~~\partial \Omega
\end{aligned}\right.
\end{equation*}
respectively, where $f_1,f_2\in L^{q}(\Omega)$ with $q>n/p$ and $1/p+1/q\leq 1$. Then
\begin{equation}\label{e.2.2}
 \sup_{\Omega} (u-v)\leq \|(g_1-g_2)^+\|_{L^{\infty}(\partial \Omega)}
 +C|\Omega|^{\frac{p}{n(p-1)}-\frac{1}{q(p-1)}}\|(f_1-f_2)^-\|^{\frac{1}{p-1}}_{L^{q}(\Omega)},
\end{equation}
where $C$ depends only on $n,p$ and $q$.
\end{lemma}
\proof For any $r\geq 0$, let $A(r)=\left\{x\in \Omega: u(x)-v(x)-\|(g_1-g_2)^+\|_{L^{\infty}(\partial \Omega)}\geq r\right\}$. Given $k\geq 0$, define
\begin{equation*}
  w(x)=\left\{\begin{aligned}
  &u(x)-v(x)-\|(g_1-g_2)^+\|_{L^{\infty}(\partial \Omega)}-k,
  &&~x\in A(k);\\
  &0 &&~x\in \Omega\backslash A(k).
  \end{aligned}\right.
\end{equation*}
Note that $w\geq 0$ in $\Omega$ and $w=0$ on $\partial \Omega$ for any $k\geq 0$. Then by the Poincar\'{e} inequality and basic calculation, we have
\begin{equation*}
  \begin{aligned}
    |A(k)|^{-p/n}\int_{A_k} |w|^p&\leq C\int_{A_k} |Dw|^p\\
    &=C\int_{A_k} |Du-Dv|^p\\
    &\leq C\int_{A_k} \langle|Du|^{p-2}Du-|Dv|^{p-2}Dv, Du-Dv\rangle\\
    &=C\int_{A_k} \langle|Du|^{p-2}Du-|Dv|^{p-2}Dv, Dw\rangle\\
    &=C\int_{A_k} (f_2-f_1) w\\
    &\leq C\left(\int_{A_k} w^p\right)^{1/p}\left(\int_{A_k} |(f_1-f_2)^-|^{p/(p-1)}\right)^{(p-1)/p}.
  \end{aligned}
\end{equation*}
Hence,
\begin{equation*}
  \begin{aligned}
\|w\|^{p-1}_{L^{p}(A(k))}
    &\leq C|A(k)|^{p/n}\|(f_1-f_2)^-\|_{L^{p/(p-1)}(A(k))}\\
    &\leq C|A(k)|^{p/n+(p-1)/p-1/q}\|(f_1-f_2)^-\|_{L^{q}(A(k))}.
  \end{aligned}
\end{equation*}
For any $h\geq k$,
\begin{equation*}
  \begin{aligned}
|A(h)|^{1/p}(h-k)&= \|h-k\|_{L^{p}(A(h))}\\
&\leq \|w\|_{L^{p}(A(h))}\\
&\leq \|w\|_{L^{p}(A(k))}\\
&\leq C|A(k)|^{p/(n(p-1))+1/p-1/(q(p-1))}\|(f_1-f_2)^-\|^{1/(p-1)}_{L^{q}(A(k))}.
  \end{aligned}
\end{equation*}
Thus,
\begin{equation*}
  |A(h)|\leq \frac{C\|(f_1-f_2)^-\|^{p/(p-1)}_{L^{q}(A(k))}}{(h-k)^p}|A(k)|^{p^2/(n(p-1))+1-p/(q(p-1))}.
\end{equation*}

Now, we apply \Cref{l-2} with $\phi(r)=|A(r)|$, $\beta=1+p^2/(n(p-1))-p/(q(p-1))>1$ and $\gamma=p$. Note that $\phi(0)\leq |\Omega|$. Then
\begin{equation*}
  \phi(r)=0~~~\mbox{for}~
  r=C2^{\frac{\beta}{\beta-1}}\|(f_1-f_2)^-\|^{1/(p-1)}_{L^{q}(\Omega)}
  |\Omega|^{\frac{p}{n(p-1)}-\frac{1}{q(p-1)}}.
\end{equation*}
That is, \cref{e.2.2} holds.  \qed~\\

Now, we can prove the boundary H\"{o}lder regularity for the $p$-Poisson equations.

\noindent \textbf{Proof of \Cref{t-5}.} Without loss of generality, we assume that $g(0)=0$. Let $M=\|u\|_{L^{\infty }(\Omega\cap B_1)}+\|f\|_{C_{q}^{-p,(p-1)\alpha}(0)}+[g]_{C^{\alpha}(0)}$ and $\Omega_{r}=\Omega\cap B_{r}$. To prove that $u$ is $C^{\alpha}$ at 0, we only need to show the following:

there exist constants $0<\delta<1/4$, $0< \eta < 1$ (depending only on $n,\alpha$ and $p$) and $\hat{C}\geq 1$ (depending only on $n,\alpha,p$ and $q$) such that for all $k\geq 0$,

\begin{equation}\label{e-ful-discrete}
\|u\|_{L^{\infty }(\Omega _{\eta ^{k}})}\leq \hat{C} M \eta^{k\alpha}.
\end{equation}

We prove \cref{e-ful-discrete} by induction. For $k=0$, it holds clearly. Suppose that it holds for $k$. We need to prove that it holds for $k+1$.

Let $r=\eta ^{k}/2$. Then there exists a new coordinate system denoted by $\{x_1,...,x_n \}$ again such that
\begin{equation*}
B_r\cap \Omega \subset B_r \cap\{x_n >-\delta r\}.
\end{equation*}
Let $\tilde{B}^{+}_{r}=B^{+}_{r}-\delta r e_n $, $\tilde{T}_r=T_r-\delta r e_n$ and $\tilde{\Omega }_{r}=\Omega \cap \tilde{B}^{+}_{r}$. Then $\Omega _{r/2}\subset \tilde{\Omega }_{r}\subset \Omega_{\eta^k}$.

Define
\begin{equation*}
v(x)=\hat{C}M\eta^{k\alpha}v_0(x/r)
\end{equation*}
and by choosing $\beta$ large enough, $v$ satisfies
\begin{equation*}
\left\{\begin{aligned}
 &\mathrm{div}(|Dv|^{p-2}Dv)\leq 0 &&\mbox{in}~~\tilde{B}^{+}_{r}; \\
 &v\geq 0  &&\mbox{in}~~\tilde{B}^{+}_{r};\\
 &v\geq \hat{C} M \eta^{k\alpha}&&\mbox{on}~~\partial \tilde{B}^{+}_{r}\backslash \tilde{T}_{r};\\
  &v(-\delta r e_n)=0.
\end{aligned}
\right.
\end{equation*}

As before, take $\delta = \eta$ and
\begin{equation}\label{e-v}
v\leq C_1\eta ^{1-\alpha }\cdot \hat{C}M\eta ^{(k+1)\alpha}~\mbox{on}~\left\{(x',x_n): x'=0,-\delta r\leq x_n\leq 2\eta r\right\},
\end{equation}
where $C_1$ depends only on $n$ and $p$. In addition, by \Cref{l-2}, we have for $x\in \tilde{\Omega }_{r}$,
\begin{equation}\label{e-w}
\begin{aligned}
u(x)-v(x) &\leq\|g\|_{L^{\infty }(\partial \Omega \cap \tilde{B}^{+}_{r})}
+Cr^{\frac{1}{p-1}(p-\frac{n}{q})}\|f\|^{\frac{1}{p-1}}_{L^q(\tilde\Omega_{r})}\\
&\leq M(2r)^{\alpha}+CM(2r)^{\alpha}=\frac{C_2}{\hat{C}\eta^{\alpha}}\cdot \hat{C}M\eta^{(k+1)\alpha},
\end{aligned}
\end{equation}
where $C_2$ depends only on $n,p$ and $q$.

Take $\eta\leq 1/4$ small enough such that
\begin{equation*}
  C_1\eta ^{1-\alpha }\leq 1/2.
\end{equation*}
Next, take $\hat{C}$ large enough such that
\begin{equation*}
  \frac{C_2}{\hat{C}\eta^{\alpha}}\leq 1/2.
\end{equation*}
Then by combining with \cref{e-v} and \cref{e-w}, we have
\begin{equation*}
 u\leq \hat{C}M\eta ^{(k+1)\alpha}~~~~\mbox{on}~\left\{(x',x_n): x'=0,-\delta r\leq x_n\leq 2\eta r\right\}.
\end{equation*}

By translating $v$ to proper positions and similar arguments, we obtain
\begin{equation*}
 \sup_{\Omega_{\eta ^{k+1}}}u\leq \hat{C}M\eta ^{(k+1)\alpha}.
\end{equation*}
The proof for
\begin{equation*}
\inf_{\Omega _{\eta^{k+1}}} u\geq -\hat{C}M\eta ^{(k+1)\alpha}
\end{equation*}
is similar and we omit here. Therefore,
\begin{equation*}
\|u\|_{L^{\infty}(\Omega _{\eta^{k+1}})}\leq \hat{C}M\eta ^{(k+1)\alpha}.
\end{equation*}
By induction, the proof is completed. \qed~\\

\bibliographystyle{model4-names}
\bibliography{PDE}

%
%
%
%
%
%
%
%

\end{document}